\documentclass[12pt]{article}

\usepackage[T1]{fontenc}
\usepackage[utf8]{inputenc}
\usepackage[margin=1in]{geometry}
\usepackage{setspace}
\usepackage{amsmath,amssymb,amsthm,amsfonts,mathrsfs}
\usepackage{enumitem}
\usepackage[colorlinks=true,linkcolor=blue,citecolor=blue,urlcolor=blue]{hyperref}

\newtheorem{theorem}{Theorem}[section]
\newtheorem{lemma}[theorem]{Lemma}
\newtheorem{proposition}[theorem]{Proposition}
\newtheorem{corollary}[theorem]{Corollary}
\newtheorem{conjecture}[theorem]{Conjecture}

\newtheorem{remark}[theorem]{Remark}
\newtheorem{example}[theorem]{Example}

\newcommand{\N}{\mathbb{N}}
\newcommand{\R}{\mathbb{R}}
\newcommand{\C}{\mathbb{C}}
\newcommand{\Z}{\mathbb{Z}}
\newcommand{\Q}{\mathbb{Q}}
\newcommand{\T}{\mathbb{T}}
\newcommand{\cS}{\mathcal{S}}

\begin{document}

\title{A Trichotomy and Rigidity Constraints for the HRT Conjecture in the Mixed-Integer Case}
\author{Vignon Oussa}
\date{\today}
\maketitle

\begin{abstract}
We consider the HRT conjecture in the mixed-integer setting, where a finite configuration in $\R^d\times\R^d$ consists of $N-1$ points in $\Z^d\times\Z^d$ and one point $(\alpha,\beta)$ outside the lattice. Assuming a linear dependence among the corresponding time-frequency shifts of a nonzero Schwartz function, we apply the Zak transform to obtain a cocycle over translation by $\gamma=(-\alpha,\beta)$ on $\T^{2d}$ and study the orbit closure
\[
H=\overline{\{n\gamma \bmod \Z^{2d}:n\in\Z\}}.
\]
We show that this reduction yields a trichotomy. The dense-orbit case is impossible because a Zak zero propagates to a dense zero set, forcing the Zak transform to vanish identically. The finite-orbit case reduces to a rational configuration, and hence to the lattice case covered by Linnell's theorem. Thus any mixed-integer counterexample for a Schwartz window must occur in the infinite proper case. For that remaining case, we prove that the nonvanishing set of the Zak transform is $H$-saturated, that the averaged logarithmic growth of the modulus cocycle along $H$ exists and vanishes identically, and that the restriction to each nonvanishing $H_0$-coset satisfies a smooth cohomological equation. This yields small-divisor compatibility conditions for the induced translation on $H_0$. We further obtain an arithmetic rigidity condition. These results isolate a collection of necessary dynamical, cohomological, and arithmetic constraints that any mixed-integer counterexample must satisfy.
\end{abstract}

\begin{center}
\textit{Dedicated to the memory of Jean-Pierre Gabardo}
\end{center}

\section{Introduction}

Gabor theory studies systems generated by translation and modulation of a
single window function. For $x,y\in\R^d$, define
\[
T_xf(t)=f(t-x),\qquad M_yf(t)=e^{-2\pi i\langle y,t\rangle}f(t).
\]
The family $\{M_yT_xg\}$ is the basic object of time-frequency analysis and
may be viewed as a discretized orbit of the Schr\"odinger representation of the
Heisenberg group; see, for example, \cite{Grochenig2001,heil2010basis}. The
subject contains many of the central results of modern harmonic analysis,
including density theorems, Wexler--Raz biorthogonality, Janssen's
representation, Ron--Shen duality, and the Balian--Low theorem
\cite{Grochenig2001}.

Jean-Pierre Gabardo made deep and lasting contributions to Gabor analysis,
frame theory, and related parts of harmonic analysis. A number of his papers
either substantially extend foundational results or reveal new connections
between frame theory, invariant subspaces, and Fourier analysis; see
\cite{Gabardo2004,Gabardo2012,Gabardo2014Measure,Gabardo2015Gabor,
Gabardo2015SelfAffine,Gabardo2018Sampling,Gabardo2018Weighted,
Gabardo2019Decomposition,Gabardo2020Frames,Gabardo2020LocalFourier,
Gabardo2020OpenSet,Gabardo2021}. Among other contributions, Gabardo and Han
\cite{Gabardo2004} developed an operator-algebraic framework for Gabor frames
restricted to closed time-frequency-invariant subspaces of $L^2(\R)$. Their
work gives a precise description of uniqueness of Weyl--Heisenberg duals in
terms of the density parameter and clarifies the interplay among completeness,
redundancy, and duality in invariant subspaces.

I first became acquainted with Jean-Pierre Gabardo through his work in Gabor
and frame theory. As an abstract harmonic analyst, what drew me especially to
his research was the unusual breadth with which he moved between abstract and
applied harmonic analysis. That breadth is visible throughout his publications
and remains an important part of his mathematical legacy
\cite{Gabardo2004,Gabardo2012,Gabardo2024}. Over the years I had the privilege
of meeting Jean-Pierre at several conferences. I remember in particular the
\href{https://my.vanderbilt.edu/iccha7/}{7th International Conference on
Computational Harmonic Analysis}, held at Vanderbilt University in Nashville in
2018. During that meeting I spent meaningful time with Jean-Pierre Gabardo and
Chun-Kit Lai. At the time I was actively thinking about discretization problems
for unitary representations of solvable Lie groups
\cite{oussa2025lieframe,Oussa2017,Oussa2018,Oussa2019,Oussa2024}. Some of the
group actions appearing there could be represented in part by pointwise
multiplication by nonlinear exponential phases, and I suspected that these
ideas might interact with Jean-Pierre's work on exponential frames. During a
conversation between sessions, Jean-Pierre immediately pushed my preliminary
observations further. Together with Chun-Kit Lai, we formulated a question of
common interest that eventually led to our paper \cite{Gabardo2021}. That
experience remains, for me, a vivid example of Jean-Pierre's generosity with
both ideas and time.

Unfortunately, the years that followed did not allow the collaboration to grow
in the way I had hoped. The COVID-19 pandemic, together with the ordinary
demands of teaching and family life, interrupted that line of work. I never had
the chance to see Jean-Pierre again after the 2018 conference, and the news of
his passing, relayed by Deguang Han in the fall of 2024, was profoundly sad.
It is therefore with genuine respect and gratitude that I dedicate this paper
on the HRT conjecture to his memory.

The HRT conjecture, formulated by Heil, Ramanathan, and Topiwala in 1996,
asserts the linear independence of every finite set of distinct time-frequency
translates of a nonzero function.

\begin{conjecture}[HRT]
Let
\[
\Lambda=\{(x_k,y_k):1\leq k\leq N\}\subset \R^d\times \R^d
\]
be a finite set of distinct points, and let $0\neq f\in L^2(\R^d)$. Then the
finite Gabor system
\[
\mathcal{G}(f,\Lambda):=\{M_{y_k}T_{x_k}f:1\leq k\leq N\}
\]
is linearly independent.
\end{conjecture}

The original paper \cite{heil1996linear} established the conjecture for several
important families of configurations, including regularly spaced collinear
points with one additional point in dimension one. A major advance came in
Linnell's work \cite{linnell1999neumann}, which proves the conjecture whenever
$\Lambda$ is contained in a full-rank lattice in $\R^{2d}$ and, in particular,
whenever $N\leq 3$. Since then, many special cases have been resolved. For
example, Demeter and Zaharescu \cite{demeter2012hrt} settled the case of
four-point configurations forming a trapezoid in the plane. There are also
results for windows with strong decay or asymptotic regularity. Bownik and
Speegle \cite{bownik2013linear} proved linear independence for broad classes of
subexponentially decaying windows, while Benedetto and Bourouihiya
\cite{benedetto2015linear} obtained strong results for windows coming from
Hardy fields, ultimately analytic functions, and several other asymptotic
classes. Okoudjou \cite{okoudjou2019extension} proposed an inductive viewpoint,
asking when a configuration for which HRT is known can be enlarged by one point
without losing linear independence.

The present paper is motivated by the mixed-integer situation studied recently
with Okoudjou \cite{okoudjou2025hrt}. There one considers a configuration with
$N-1$ points in $\Z^d\times\Z^d$ and one point outside the lattice. In
dimension one, this includes configurations such as
\[
\{(0,0),(1,0),(0,1),(\sqrt2,\sqrt2)\}
\quad\text{and}\quad
\{(0,0),(1,0),(0,1),(\sqrt2,\sqrt3)\}.
\]
Applying the Zak transform produces a cocycle over translation on the torus, and
the geometry of the orbit closure leads naturally to a trichotomy: the orbit
may be dense, finite, or infinite and proper. The first two cases are favorable
for HRT. The purpose of this paper is to isolate the precise constraints that
remain in the infinite proper case. Besides vanishing averaged logarithmic
growth and arithmetic rigidity on the identity component, a Fourier analysis of
the modulus cocycle on $H_0$-cosets yields a fiberwise cohomological equation
with small divisors. The zero Fourier mode remains free, so the modulus
cocycle alone does not force a contradiction at a zero coset. What it does
provide is a sharp compatibility condition between the Fourier spectrum of the
fiber cocycle and the arithmetic of the translation direction. In particular,
this mechanism can exclude exponentially Liouville orbit directions once the
fiber cocycle retains sufficient Fourier mass along the resonant modes.

\section{Zak-transform reduction and a uniform ergodic average}

Throughout the paper we write
\[
\T^n:=\R^n/\Z^n
\]
for the $n$-torus. We freely identify periodic functions on $\R^n$ with
functions on $\T^n$.

Let
\[
\Lambda=
\Bigl\{(x_k,y_k)\in\Z^d\times\Z^d:1\le k\le N-1\Bigr\}
\cup\{(\alpha,\beta)\},
\qquad
(\alpha,\beta)\notin \Z^d\times\Z^d.
\]
Assume that there exist $0\neq f\in\cS(\R^d)$ and coefficients
$c_1,\dots,c_{N-1}\in\C$ such that
\begin{equation}
\label{eq:dependence}
\sum_{k=1}^{N-1} c_k M_{y_k}T_{x_k}f=M_\beta T_\alpha f.
\end{equation}
We define the Zak transform by
\[
Zf(t,\omega):=\sum_{\kappa\in\Z^d}e^{-2\pi i\langle \omega,\kappa\rangle}
f(t+\kappa),
\qquad (t,\omega)\in\R^d\times\R^d.
\]
For Schwartz functions, this series converges absolutely and defines a smooth
function on $\R^{2d}$. It satisfies the quasi-periodicity relations
\begin{equation}
\label{eq:zak-quasi}
Zf(t+k,\omega)=e^{2\pi i\langle \omega,k\rangle}Zf(t,\omega),
\qquad
Zf(t,\omega+\ell)=Zf(t,\omega),
\end{equation}
for every $k,\ell\in\Z^d$; see \cite[Chapter~8]{Grochenig2001}. In
particular, $|Zf|$ is $\Z^{2d}$-periodic and therefore descends to a continuous
function on $\T^{2d}$.

\begin{lemma}[Zak-transform cocycle]
\label{lem:zak-cocycle}
Let
\[
p(t,\omega):=\sum_{k=1}^{N-1}c_k e^{-2\pi i\langle y_k,t\rangle}
 e^{-2\pi i\langle \omega,x_k\rangle},
\qquad
\gamma:=(-\alpha,\beta).
\]
Then for every $(t,\omega)\in\R^{2d}$ one has
\begin{equation}
\label{eq:phase-cocycle}
Zf(t-\alpha,\omega+\beta)=e^{2\pi i\langle t,\beta\rangle}p(t,\omega)Zf(t,\omega).
\end{equation}
Consequently, if
\[
F:=|Zf|,
\qquad
q:=|p|,
\]
then
\begin{equation}
\label{eq:modulus-cocycle}
F(z+\gamma)=q(z)F(z)
\qquad (z\in\T^{2d}).
\end{equation}
Moreover,
\begin{equation}
\label{eq:iterated-cocycle}
F(z+n\gamma)=\Bigl(\prod_{j=0}^{n-1}q(z+j\gamma)\Bigr)F(z)
\qquad (n\in\N).
\end{equation}
\end{lemma}

\begin{proof}
Applying the Zak transform to the left-hand side of \eqref{eq:dependence} and
using that $x_k,y_k\in\Z^d$, we obtain
\[
Z(M_{y_k}T_{x_k}f)(t,\omega)
=e^{-2\pi i\langle y_k,t\rangle}e^{-2\pi i\langle \omega,x_k\rangle}Zf(t,\omega).
\]
Hence
\[
Z\Bigl(\sum_{k=1}^{N-1}c_kM_{y_k}T_{x_k}f\Bigr)(t,\omega)=p(t,\omega)Zf(t,\omega).
\]
For the right-hand side we compute directly:
\[
\begin{aligned}
Z(M_\beta T_\alpha f)(t,\omega)
&=\sum_{\kappa\in\Z^d}e^{-2\pi i\langle \omega,\kappa\rangle}
 e^{-2\pi i\langle \beta,t+\kappa\rangle}f(t+\kappa-\alpha)\\
&=e^{-2\pi i\langle t,\beta\rangle}
 \sum_{\kappa\in\Z^d}e^{-2\pi i\langle \omega+\beta,\kappa\rangle}f((t-\alpha)+\kappa)\\
&=e^{-2\pi i\langle t,\beta\rangle}Zf(t-\alpha,\omega+\beta).
\end{aligned}
\]
Combining the two identities gives \eqref{eq:phase-cocycle}. Taking absolute
values yields \eqref{eq:modulus-cocycle}, and iterating the latter gives
\eqref{eq:iterated-cocycle}.
\end{proof}

We now record two standard facts that will be used repeatedly.

\begin{lemma}
\label{lem:continuous-zak-zero}
If $0\neq f\in \cS(\R^d)$, then $Zf$ is continuous and has a zero on every
fundamental domain of $\T^{2d}$.
\end{lemma}

\begin{proof}
Continuity is immediate for Schwartz functions. The existence of a zero is a
standard property of continuous Zak transforms; see \cite[Lemma~8.4.2]{Grochenig2001}.
\end{proof}
The following is a standard ergodic theorem
\begin{lemma}\label{lem:uniform-averages}
Let $H$ be a compact abelian group, let $\gamma\in H$, and assume that the
cyclic subgroup generated by $\gamma$ is dense in $H$. Then for every
$\varphi\in C(H)$,
\[
\frac1n\sum_{j=0}^{n-1}\varphi(h+j\gamma)
\longrightarrow
\int_H\varphi\,dm_H
\qquad\text{uniformly in }h\in H,
\]
where $m_H$ denotes normalized Haar measure on $H$.
\end{lemma}

\section{The orbit trichotomy}

Define
\begin{equation}
\label{eq:def-H}
H:=\overline{\{n\gamma\bmod \Z^{2d}:n\in\Z\}}\subset \T^{2d}.
\end{equation}
Since $H$ is a closed subgroup of the torus, it is a compact abelian Lie group;
see, for example, \cite[Chapter~2]{HofmannMorris2013}.

\begin{proposition}[Orbit saturation of the nonvanishing set]
\label{prop:orbit-saturation}
For every $z\in S$, the entire coset $z+H$ is contained in $S$. Moreover,
\[
p(w)\neq 0\qquad\text{for every }w\in z+H.
\]
\end{proposition}

\begin{proof}
Let $Z:=\T^{2d}\setminus S=\{F=0\}$. By \eqref{eq:modulus-cocycle}, the zero
set $Z$ is forward invariant under translation by $\gamma$: if $u\in Z$, then
$u+n\gamma\in Z$ for every $n\geq 0$. Fix $z\in S$ and suppose first that there exists $w\in z+H$ with $F(w)=0$. The closure of the forward orbit
$\{w+n\gamma:n\geq 0\}$ is exactly $w+H=z+H$. Since $Z$ is closed and contains
that forward orbit, it contains all of $z+H$, contradicting $z\in S$. Thus, $F$ is strictly positive on $z+H$, so $z+H\subset S$. Now suppose that $p(w_0)=0$ for some $w_0\in z+H$. Then
\eqref{eq:modulus-cocycle} gives
\[
F(w_0+\gamma)=q(w_0)F(w_0)=0.
\]
Since $w_0+\gamma\in z+H$, the first part of the proof gives a contradiction.
Therefore $p$ has no zero on $z+H$.
\end{proof}

The previous proposition allows us to define logarithmic averages on every
$H$-coset contained in $S$.

\begin{proposition}[Averaged logarithmic growth]
\label{prop:theta}
For each $z\in S$, define
\begin{equation}
\label{eq:def-theta}
\Theta(z):=\int_H\log|p(z+h)|\,dm_H(h).
\end{equation}
Then the following hold.
\begin{enumerate}[label=\textup{(\roman*)}]
\item The function $\Theta$ is well defined on $S$ and is constant on every
coset $z+H$.
\item For every $z\in S$,
\begin{equation}
\label{eq:theta-limit}
\lim_{n\to\infty}\frac{\log F(z+n\gamma)-\log F(z)}{n}=\Theta(z).
\end{equation}
\item The function $\Theta$ is real-analytic on each connected component of
$S$.
\item One has
\begin{equation}
\label{eq:theta-zero}
\Theta(z)=0\qquad\text{for every }z\in S.
\end{equation}
\end{enumerate}
\end{proposition}

\begin{proof}
Fix $z\in S$. By Proposition~\ref{prop:orbit-saturation}, the function
$h\mapsto p(z+h)$ is nonvanishing on the compact group $H$, so
$h\mapsto \log|p(z+h)|$ is continuous. This proves that $\Theta(z)$ is well
defined. Haar invariance immediately gives
\[
\Theta(z+h_0)=\Theta(z)\qquad (h_0\in H),
\]
so $\Theta$ is constant on each $H$-coset. To prove \eqref{eq:theta-limit}, observe from \eqref{eq:iterated-cocycle} that
\[
\frac{\log F(z+n\gamma)-\log F(z)}{n}
=\frac1n\sum_{j=0}^{n-1}\log q(z+j\gamma).
\]
Apply Lemma~\ref{lem:uniform-averages} to the compact group $H$ and to the
continuous function
\[
\varphi_z(h):=\log q(z+h)=\log|p(z+h)|.
\]
Because the subgroup generated by $\gamma$ is dense in $H$ by definition,
Lemma~\ref{lem:uniform-averages} yields
\[
\frac1n\sum_{j=0}^{n-1}\log q(z+j\gamma)
\longrightarrow
\int_H\log q(z+h)\,dm_H(h)=\Theta(z),
\]
which proves \eqref{eq:theta-limit}.

We next prove real-analyticity. Define
\[
g(z):=\min_{h\in H}|p(z+h)|.
\]
Since $p$ is continuous and $H$ is compact, $g$ is continuous on $\T^{2d}$. By
Proposition~\ref{prop:orbit-saturation}, $g(z)>0$ for every $z\in S$. Hence, for
each $z_0\in S$, there exists a neighborhood $U$ of $z_0$ and a constant
$\delta>0$ such that
\[
|p(z+h)|\geq \delta\qquad (z\in U,\;h\in H).
\]
The function $(z,h)\mapsto |p(z+h)|^2$ is real-analytic on $U\times H$ and
strictly positive there, so $(z,h)\mapsto \log|p(z+h)|$ is real-analytic on
$U\times H$. Integrating over the compact parameter space $H$ shows that
$\Theta$ is real-analytic on $U$. Since $z_0$ was arbitrary, $\Theta$ is
real-analytic on each connected component of $S$. Finally, let $z\in S$. If $\Theta(z)>0$, then by \eqref{eq:theta-limit} there
exists $N$ such that for all $n\geq N$,
\[
\log F(z+n\gamma)-\log F(z)\geq \frac{\Theta(z)}{2}n.
\]
Hence
\[
F(z+n\gamma)\geq F(z)e^{\Theta(z)n/2}\to\infty,
\]
which is impossible because $F$ is continuous on the compact torus $\T^{2d}$.
If $\Theta(z)<0$, choose a sequence $n_j\to\infty$ with $n_j\gamma\to 0.$ Then continuity gives
\[
F(z+n_j\gamma)\to F(z)>0.
\]
On the other hand, \eqref{eq:theta-limit} implies that for all sufficiently
large $j$,
\[
\log F(z+n_j\gamma)-\log F(z)\leq \frac{\Theta(z)}{2}n_j,
\]
so
\[
F(z+n_j\gamma)\leq F(z)e^{\Theta(z)n_j/2}\to 0,
\]
a contradiction. Therefore neither inequality is possible, and
$\Theta(z)=0$.
\end{proof}

We can now state the trichotomy in its final form.

\begin{theorem}[Trichotomy for mixed-integer configurations]
\label{thm:trichotomy}
Assume \eqref{eq:dependence} with $0\neq f\in\cS(\R^d)$, and let $H$ be given by
\eqref{eq:def-H}. Then exactly one of the following three alternatives occurs.
\begin{enumerate}[label=\textup{(\roman*)}]
\item \textbf{Dense-orbit case:} $H=\T^{2d}$. This case is impossible.
\item \textbf{Finite-orbit case:} $H$ is finite. In this case
$(\alpha,\beta)\in\Q^{2d}\setminus\Z^{2d}$, the set $\Lambda$ is contained in a
full-rank lattice in $\R^{2d}$, and therefore Linnell's theorem implies
$f=0$.
\item \textbf{Infinite proper case:} $H$ is infinite and proper. In this case
$S$ is a nonempty open $H$-invariant subset of $\T^{2d}$, and the function
$\Theta$ defined by \eqref{eq:def-theta} is well defined on $S$ and satisfies
\eqref{eq:theta-limit} and \eqref{eq:theta-zero}.
\end{enumerate}
In particular, if a nontrivial mixed-integer counterexample to the HRT
conjecture exists for a Schwartz window, then it must belong to the infinite
proper case and satisfy $\Theta\equiv 0$ on $S$.
\end{theorem}

\begin{proof}
The three alternatives are mutually exclusive and exhaustive because $H$ is a
closed subgroup of $\T^{2d}$. Suppose first that $H=\T^{2d}$. By
Lemma~\ref{lem:continuous-zak-zero}, the continuous function $Zf$ has a zero,
say at $z_0\in\T^{2d}$. Then the forward orbit $\{z_0+n\gamma:n\geq 0\}$ is
dense in $\T^{2d}$, and \eqref{eq:phase-cocycle} shows that the zero set of
$Zf$ is forward invariant under $\gamma$. Hence the zero set of $Zf$ is dense.
Since the zero set is also closed, it follows that $Zf\equiv 0$, hence $f=0$,
a contradiction. Therefore the dense-orbit case is impossible. Suppose next that $H$ is finite. Then the class of $\gamma$ in $\T^{2d}$ is
torsion, so there exists $m\in\N$ such that $m\gamma\in\Z^{2d}$. Equivalently,
$\gamma\in\Q^{2d}$, that is, $(\alpha,\beta)\in\Q^{2d}$. Choose
$M\in\N$ such that $M\alpha\in\Z^d$ and $M\beta\in\Z^d$. Then
\[
\Lambda\subset M^{-1}\Z^d\times M^{-1}\Z^d,
\]
which is a full-rank lattice in $\R^{2d}$. Linnell's theorem
\cite{linnell1999neumann} now implies that $f=0$.

The third case is exactly Proposition~\ref{prop:theta}.
\end{proof}

\section{The identity component and the $m$-step cocycle}

We now assume throughout that $H$ is infinite and proper. Let $H_0$ denote the
identity component of $H$. Since $H$ is a closed subgroup of a torus, it is a
compact abelian Lie group. Hence $H_0$ is a torus and the quotient $H/H_0$ is
finite \cite[Chapter~2]{HofmannMorris2013}. Set
\[
m:=|H/H_0|.
\]
Because $H/H_0$ is generated by the class of $\gamma$, it is finite cyclic.
Therefore
\[
m\gamma\in H_0.
\]
The next lemma identifies the connected component as the orbit closure of
$m\gamma$.

\begin{lemma}
\label{lem:H0-generated}
Let
\[
K:=\overline{\{nm\gamma:n\in\Z\}}\subset H.
\]
Then $K=H_0$.
\end{lemma}

\begin{proof}
Since $m\gamma\in H_0$ and $H_0$ is closed, we have $K\subset H_0$. On the
other hand, the quotient $H/K$ is generated by the class of $\gamma$, and this
class has order dividing $m$ because $m\gamma\in K$. Hence $H/K$ is a finite
cyclic group. Since $H_0/K$ is a connected subgroup of the finite group $H/K$,
it must be trivial. Thus $H_0=K$.
\end{proof}

Fix $\lambda\in S$. By Proposition~\ref{prop:orbit-saturation}, the entire coset
$\lambda+H$ lies in $S$ and $p$ is nowhere zero on $\lambda+H$. In particular,
if
\[
C_\lambda:=\lambda+H_0,
\]
then both $Zf$ and $p$ are nowhere zero on $C_\lambda$. Define
\[
P_m(z):=\prod_{j=0}^{m-1}p(z+j\gamma),
\qquad
Q_m(z):=|P_m(z)|=\prod_{j=0}^{m-1}q(z+j\gamma).
\]
Then $P_m$ is smooth and nowhere zero on $\lambda+H$, hence on $C_\lambda$.
Iterating \eqref{eq:phase-cocycle} $m$ times gives
\begin{equation}
\label{eq:m-step-cocycle}
Zf(z+m\gamma)=E_m(z)P_m(z)Zf(z),
\end{equation}
where, for $z=(t,\omega)$,
\begin{equation}
\label{eq:def-Em}
E_m(z):=\exp\!\left(2\pi i\left[m\langle t,\beta\rangle-
\frac{m(m-1)}{2}\langle \alpha,\beta\rangle\right]\right).
\end{equation}
Taking absolute values yields
\begin{equation}
\label{eq:m-step-modulus}
F(z+m\gamma)=Q_m(z)F(z).
\end{equation}

\section{Fiberwise equation and small divisors}

This section records the precise Fourier-theoretic information carried by the
modulus equation in the infinite proper case.

\begin{proposition}[Fiberwise equation]
\label{prop:cohomological}
Assume that $H$ is infinite and proper. Let $r=\dim H_0$, and choose a Lie
group isomorphism
\[
\Phi:H_0\longrightarrow \T^r.
\]
Set
\[
\omega:=\Phi(m\gamma)\in\T^r.
\]
Then $\omega$ has dense orbit in $\T^r$. Fix $\lambda\in S$ and define
\[
u_\lambda(x):=\log F(\lambda+\Phi^{-1}(x)),
\qquad
b_\lambda(x):=\log Q_m(\lambda+\Phi^{-1}(x)),
\qquad x\in\T^r.
\]
Then:
\begin{enumerate}[label=\textup{(\roman*)}]
\item $u_\lambda\in C^\infty(\T^r)$ and $b_\lambda$ is real-analytic on
$\T^r$;
\item one has the cohomological equation
\begin{equation}
\label{eq:cohomological}
u_\lambda(x+\omega)-u_\lambda(x)=b_\lambda(x)
\qquad (x\in\T^r);
\end{equation}
\item the mean of $b_\lambda$ vanishes:
\begin{equation}
\label{eq:b-mean-zero}
\widehat b_\lambda(0)=\int_{\T^r}b_\lambda(x)\,dx=0;
\end{equation}
\item for every $k\in\Z^r\setminus\{0\}$,
\begin{equation}
\label{eq:uk-formula}
\widehat u_\lambda(k)=\frac{\widehat b_\lambda(k)}{e^{2\pi i\langle k,\omega\rangle}-1};
\end{equation}
\item there exists a constant $c(\lambda)\in\R$ such that
\begin{equation}
\label{eq:u-fourier-series}
u_\lambda(x)=c(\lambda)+
\sum_{k\in\Z^r\setminus\{0\}}
\frac{\widehat b_\lambda(k)}{e^{2\pi i\langle k,\omega\rangle}-1}
 e^{2\pi i\langle k,x\rangle},
\end{equation}
with convergence in $C^\infty(\T^r)$.
\end{enumerate}
\end{proposition}

\begin{proof}
Since $C_\lambda\subset S$, the function $F=|Zf|$ is smooth and strictly
positive on $C_\lambda$, hence $u_\lambda\in C^\infty(\T^r)$. The function
$Q_m$ is real-analytic on $C_\lambda$ and strictly positive there, so
$b_\lambda$ is real-analytic. The density of the cyclic subgroup generated by $\omega$ follows from
Lemma~\ref{lem:H0-generated}. Indeed, the dense subgroup generated by $m\gamma$
in $H_0$ is carried by $\Phi$ to the cyclic subgroup generated by $\omega$ in
$\T^r$. The $m$-step modulus equation \eqref{eq:m-step-modulus} gives
\[
F(\lambda+\Phi^{-1}(x+\omega))=Q_m(\lambda+\Phi^{-1}(x))F(\lambda+\Phi^{-1}(x)),
\]
which is exactly \eqref{eq:cohomological} after taking logarithms.
Integrating \eqref{eq:cohomological} over $\T^r$ yields
\eqref{eq:b-mean-zero}. Now take Fourier coefficients in \eqref{eq:cohomological}. For each
$k\in\Z^r$,
\[
\widehat{u_\lambda(\cdot+\omega)}(k)=e^{2\pi i\langle k,\omega\rangle}\widehat u_\lambda(k),
\]
so
\[
\bigl(e^{2\pi i\langle k,\omega\rangle}-1\bigr)\widehat u_\lambda(k)=\widehat b_\lambda(k).
\]
If $k\neq 0$, then $e^{2\pi i\langle k,\omega\rangle}\neq 1$, because otherwise
the corresponding character would vanish on the dense cyclic subgroup generated
by $\omega$, hence on all of $\T^r$. This proves \eqref{eq:uk-formula}. The
case $k=0$ is exactly \eqref{eq:b-mean-zero}. Finally,
\eqref{eq:u-fourier-series} is just the Fourier expansion of the smooth
function $u_\lambda$, with
\[
c(\lambda)=\widehat u_\lambda(0)=\int_{\T^r}u_\lambda(x)\,dx.
\]
\end{proof}

\begin{proposition}[Small-divisor compatibility]
\label{prop:small-divisor}
In the setting of Proposition~\ref{prop:cohomological}, for every integer
$N\geq 0$ one has
\begin{equation}
\label{eq:small-divisor-compatibility}
\sup_{k\in\Z^r\setminus\{0\}}(1+|k|)^N
\frac{|\widehat b_\lambda(k)|}{|e^{2\pi i\langle k,\omega\rangle}-1|}<\infty.
\end{equation}
Equivalently,
\begin{equation}
\label{eq:small-divisor-compatibility-distance}
\sup_{k\in\Z^r\setminus\{0\}}(1+|k|)^N
\frac{|\widehat b_\lambda(k)|}{\|\langle k,\omega\rangle\|_{\T}}<\infty,
\end{equation}
where $\|t\|_{\T}:=\inf_{n\in\Z}|t-n|$.
In particular, if there exist $N\geq 0$ and a sequence
$k_n\in\Z^r\setminus\{0\}$ such that
\[
(1+|k_n|)^N
\frac{|\widehat b_\lambda(k_n)|}{\|\langle k_n,\omega\rangle\|_{\T}}
\longrightarrow\infty,
\]
then no nonzero Schwartz window can satisfy \eqref{eq:dependence}.
\end{proposition}

\begin{proof}
Since $u_\lambda\in C^\infty(\T^r)$, its Fourier coefficients decay rapidly:
for every $N\geq 0$,
\[
\sup_{k\in\Z^r}(1+|k|)^N|\widehat u_\lambda(k)|<\infty.
\]
Combining this with \eqref{eq:uk-formula} yields
\eqref{eq:small-divisor-compatibility}. The equivalence with
\eqref{eq:small-divisor-compatibility-distance} follows from the standard
estimate
\[
|e^{2\pi i t}-1|\asymp \|t\|_{\T}.
\]
The final claim is immediate.
\end{proof}

\section{Arithmetic rigidity on the identity component}

We continue in the infinite proper case. Fix $\lambda\in S$ and choose a
representative $\widetilde\lambda\in\R^{2d}$ of $\lambda$. Let
\[
\mathfrak h:=T_0H_0\subset\R^{2d},
\qquad
\widetilde C_\lambda:=\widetilde\lambda+\mathfrak h\subset\R^{2d}.
\]
The quotient map $\R^{2d}\to\T^{2d}$ identifies
$\widetilde C_\lambda/(\mathfrak h\cap\Z^{2d})$ with
$C_\lambda=\lambda+H_0$. By Proposition~\ref{prop:orbit-saturation}, the zero set of $Zf$ does not meet
$C_\lambda$. Equivalently, $Zf$ is nowhere zero on $\widetilde C_\lambda$.
Hence the logarithmic derivatives $L_v(Zf)$ are well defined on
$\widetilde C_\lambda$. For a smooth nowhere-zero complex-valued function $G$ and a vector
$v=(v_t,v_\omega)\in\R^d\times\R^d$, write
\[
L_vG:=\frac{D_vG}{G},
\qquad
D_v:=\sum_{j=1}^d v_{t,j}\frac{\partial}{\partial t_j}+
\sum_{j=1}^d v_{\omega,j}\frac{\partial}{\partial \omega_j}.
\]
Since $Zf$ is smooth and nowhere zero on $C_\lambda$, the logarithmic
derivatives $L_v(Zf)$ are well defined on $C_\lambda$. Applying $L_v$ to \eqref{eq:m-step-cocycle} gives
\begin{equation}
\label{eq:Lv-basic}
L_v(Zf)(z+m\gamma)-L_v(Zf)(z)=L_v(E_m)(z)+L_v(P_m)(z).
\end{equation}
Because $E_m$ depends only on the $t$-variable, \eqref{eq:def-Em} yields
\begin{equation}
\label{eq:Lv-Em}
L_v(E_m)(z)=2\pi i\,m\langle v_t,\beta\rangle.
\end{equation}
Hence
\begin{equation}
\label{eq:Lv-iterate}
L_v(Zf)(z+m\gamma)-L_v(Zf)(z)=2\pi i\,m\langle v_t,\beta\rangle+L_v(P_m)(z).
\end{equation}

The function $L_v(Zf)$ is not periodic because $Zf$ is only quasi-periodic.
Indeed, \eqref{eq:zak-quasi} implies
\[
Zf(t+k,\omega+\ell)=e^{2\pi i\langle \omega,k\rangle}Zf(t,\omega),
\qquad (k,\ell)\in\Z^d\times\Z^d,
\]
and differentiation yields
\begin{equation}
\label{eq:Lv-quasi}
L_v(Zf)(t+k,\omega+\ell)=L_v(Zf)(t,\omega)+2\pi i\langle v_\omega,k\rangle.
\end{equation}
Therefore the corrected function
\[
M_v(t,\omega):=L_v(Zf)(t,\omega)-2\pi i\langle v_\omega,t\rangle
\]
is $\Z^{2d}$-periodic. Since $z+m\gamma=(t-m\alpha,\omega+m\beta)$, we obtain
\[
\begin{aligned}
M_v(z+m\gamma)-M_v(z)
&=L_v(Zf)(z+m\gamma)-L_v(Zf)(z)
+2\pi i\,m\langle v_\omega,\alpha\rangle.
\end{aligned}
\]
Combining this with \eqref{eq:Lv-iterate}, we arrive at
\begin{equation}
\label{eq:Mv-iterate}
M_v(z+m\gamma)-M_v(z)
=2\pi i\,m\bigl(\langle v_t,\beta\rangle+\langle v_\omega,\alpha\rangle\bigr)
+L_v(P_m)(z).
\end{equation}
Integrate \eqref{eq:Mv-iterate} over the compact group $H_0$. Writing
$z=\lambda+\eta$ with $\eta\in H_0$ and denoting normalized Haar measure by
$d\mu_{H_0}$, we obtain
\begin{equation}
\label{eq:averaged-identity}
0=
2\pi i\,m\bigl(\langle v_t,\beta\rangle+\langle v_\omega,\alpha\rangle\bigr)
+\int_{H_0}L_v(P_m)(\lambda+\eta)\,d\mu_{H_0}(\eta),
\end{equation}
because $m\gamma\in H_0$ and Haar measure is translation invariant. Since $H_0$ is a connected closed subgroup of the torus, one has
\[
H_0\cong \mathfrak h/(\mathfrak h\cap\Z^{2d}),
\]
and $\mathfrak h\cap\Z^{2d}$ is a full lattice in $\mathfrak h$
\cite[Chapter~2]{HofmannMorris2013}. In particular,
$\mathfrak h\cap\Z^{2d}\neq\{0\}$ because $H$ is infinite and therefore
$H_0$ is nontrivial.

\begin{lemma}[Winding number along lattice loops]
\label{lem:winding}
Fix $a=(a_t,a_\omega)\in \mathfrak h\cap\Z^{2d}$ with $a\neq 0$. For each
$\eta\in H_0$, define
\[
c_{a,\eta}(s):=\lambda+\eta+sa\pmod{\Z^{2d}},
\qquad 0\le s\le 1,
\]
and set
\[
\kappa(a,\eta):=\frac{1}{2\pi i}\int_0^1L_a(P_m)(\lambda+\eta+sa)\,ds.
\]
Then:
\begin{enumerate}[label=\textup{(\roman*)}]
\item $\kappa(a,\eta)\in\Z$ for every $\eta\in H_0$;
\item the map $\eta\mapsto \kappa(a,\eta)$ is continuous on $H_0$ and hence
constant, say equal to $\kappa_a\in\Z$;
\item
\begin{equation}
\label{eq:kappa-identity}
\int_{H_0}L_a(P_m)(\lambda+\eta)\,d\mu_{H_0}(\eta)=2\pi i\,\kappa_a.
\end{equation}
\end{enumerate}
\end{lemma}

\begin{proof}
Because $a\in\mathfrak h$, the path $s\mapsto sa\pmod{\Z^{2d}}$ lies in $H_0$,
so $c_{a,\eta}(s)\in C_\lambda$ for all $s$. Because $a\in\Z^{2d}$, the path
is closed:
\[
c_{a,\eta}(1)=\lambda+\eta+a\equiv \lambda+\eta=c_{a,\eta}(0)
\pmod{\Z^{2d}}.
\]
Now define the closed curve
\[
\Gamma_\eta(s):=P_m(c_{a,\eta}(s))\in \C^\times.
\]
Since $P_m$ is nowhere zero on $C_\lambda$, the curve $\Gamma_\eta$ avoids the
origin. By the chain rule,
\[
\Gamma_\eta'(s)=D_aP_m(\lambda+\eta+sa),
\]
so
\[
\kappa(a,\eta)=\frac{1}{2\pi i}\int_0^1\frac{\Gamma_\eta'(s)}{\Gamma_\eta(s)}\,ds.
\]
Choose a continuous argument $\theta_\eta$ on $[0,1]$ and write
$\Gamma_\eta(s)=r_\eta(s)e^{i\theta_\eta(s)}$ with $r_\eta(s)>0$. Then
\[
\frac{\Gamma_\eta'(s)}{\Gamma_\eta(s)}=
\frac{r_\eta'(s)}{r_\eta(s)}+i\theta_\eta'(s)
\]
almost everywhere, and integration gives
\[
\int_0^1\frac{\Gamma_\eta'(s)}{\Gamma_\eta(s)}\,ds
=\log r_\eta(1)-\log r_\eta(0)+i\bigl(\theta_\eta(1)-\theta_\eta(0)\bigr).
\]
Because $\Gamma_\eta$ is closed, $r_\eta(1)=r_\eta(0)$ and
$\theta_\eta(1)-\theta_\eta(0)=2\pi n_\eta$ for some $n_\eta\in\Z$. Thus
\[
\kappa(a,\eta)=n_\eta\in\Z,
\]
which proves (i). Since $C_\lambda$ is compact and $P_m$ is nowhere zero there, there exists
$\delta>0$ such that $|P_m|\geq \delta$ on $C_\lambda$. Therefore the function
\[
(\eta,s)\longmapsto L_a(P_m)(\lambda+\eta+sa)
\]
is continuous on $H_0\times[0,1]$, and hence the map
$\eta\mapsto\kappa(a,\eta)$ is continuous. By (i), it is integer-valued on the
connected set $H_0$, so it must be constant; this proves (ii). For (iii), average the identity defining $\kappa(a,\eta)$ over $H_0$ and apply
Fubini's theorem:
\[
\begin{aligned}
\kappa_a
&=\frac{1}{2\pi i}\int_{H_0}\int_0^1L_a(P_m)(\lambda+\eta+sa)\,ds\,d\mu_{H_0}(\eta).
\end{aligned}
\]
Because $sa\in H_0$ for every $s\in[0,1]$, translation invariance of Haar
measure gives
\[
\int_{H_0}L_a(P_m)(\lambda+\eta+sa)\,d\mu_{H_0}(\eta)
=\int_{H_0}L_a(P_m)(\lambda+\eta)\,d\mu_{H_0}(\eta).
\]
Thus
\[
\kappa_a
=\frac{1}{2\pi i}\int_0^1
\left(\int_{H_0}L_a(P_m)(\lambda+\eta)\,d\mu_{H_0}(\eta)\right)ds,
\]
which is exactly \eqref{eq:kappa-identity}.
\end{proof}

\begin{proposition}[Arithmetic rigidity on the identity component]
\label{prop:arithmetic}
Assume that $H$ is infinite and proper. Let $H_0$ be the identity component of
$H$, let $m=|H/H_0|$, and let $\mathfrak h=T_0H_0$. Then for every
$a=(a_t,a_\omega)\in \mathfrak h\cap\Z^{2d}$ one has
\begin{equation}
\label{eq:arithmetic-rigidity}
m\bigl(\langle a_t,\beta\rangle+\langle a_\omega,\alpha\rangle\bigr)\in\Z.
\end{equation}
\end{proposition}

\begin{proof}
If $a=0$, the claim is trivial. Let $a\in \mathfrak h\cap\Z^{2d}$ be nonzero.
Substitute $v=a$ into \eqref{eq:averaged-identity}. Then
\[
0=2\pi i\,m\bigl(\langle a_t,\beta\rangle+\langle a_\omega,\alpha\rangle\bigr)
+\int_{H_0}L_a(P_m)(\lambda+\eta)\,d\mu_{H_0}(\eta).
\]
By Lemma~\ref{lem:winding}, the integral equals $2\pi i\,\kappa_a$ for some
$\kappa_a\in\Z$. Hence
\[
0=2\pi i\left[m\bigl(\langle a_t,\beta\rangle+\langle a_\omega,\alpha\rangle\bigr)
+\kappa_a\right],
\]
so
\[
m\bigl(\langle a_t,\beta\rangle+\langle a_\omega,\alpha\rangle\bigr)
=-\kappa_a\in\Z.
\]
This is exactly \eqref{eq:arithmetic-rigidity}.
\end{proof}

\section{Consequences and examples}

Combining Theorem~\ref{thm:trichotomy},
Proposition~\ref{prop:small-divisor}, and
Proposition~\ref{prop:arithmetic}, we obtain the main reduction statement.

\begin{corollary}
\label{cor:main-reduction}
Let
\[
\Lambda=
\Bigl\{(x_k,y_k)\in\Z^d\times\Z^d:1\le k\le N-1\Bigr\}
\cup\{(\alpha,\beta)\},
\qquad
(\alpha,\beta)\notin\Z^d\times\Z^d.
\]
If there exists a nonzero Schwartz function $f$ for which the family
$\{M_yT_xf:(x,y)\in\Lambda\}$ is linearly dependent, then the following are
necessary:
\begin{enumerate}[label=\textup{(\roman*)}]
\item the orbit closure
$H=\overline{\{n(-\alpha,\beta)\bmod\Z^{2d}:n\in\Z\}}$ is infinite and
proper;
\item with $S=\{z\in\T^{2d}:Zf(z)\neq 0\}$, the averaged logarithmic growth
function
\[
\Theta(z)=\int_H\log|p(z+h)|\,dm_H(h)
\]
is well defined on $S$ and vanishes identically there;
\item if $H_0$ denotes the identity component of $H$, if $m=|H/H_0|$, and if
$H_0\cong\T^r$ carries $m\gamma$ to $\omega\in\T^r$, then for every
$\lambda\in S$ and every $N\ge 0$ one has
\[
\sup_{k\in\Z^r\setminus\{0\}}(1+|k|)^N
\frac{|\widehat b_\lambda(k)|}{\|\langle k,\omega\rangle\|_{\T}}<\infty,
\]
where $b_\lambda$ is defined in Proposition~\ref{prop:cohomological};
\item if $H_0$ denotes the identity component of $H$ and $m=|H/H_0|$, then
\[
m\bigl(\langle a_t,\beta\rangle+\langle a_\omega,\alpha\rangle\bigr)\in\Z
\qquad
\text{for every }a=(a_t,a_\omega)\in T_0H_0\cap\Z^{2d}.
\]
\end{enumerate}
\end{corollary}

The next two concrete mixed-integer configurations are ruled out by
Proposition~\ref{prop:arithmetic} because the integrality condition
\eqref{eq:arithmetic-rigidity} fails for an explicit lattice vector in
$T_0H_0\cap\Z^{2d}$.

\begin{example}
Let
\[
\Lambda=\bigl\{((0,0),(0,0)),\;((1,0),(0,0)),\;((0,1),(0,0)),\;((0,0),(1,0)),\;(\alpha,\beta)\bigr\}
\subset\R^2\times\R^2,
\]
with $\alpha=(\sqrt2,\sqrt3)$ and $\beta=(\sqrt2,\sqrt5)$. Then
$\gamma=(-\alpha,\beta)=(-\sqrt2,-\sqrt3,\sqrt2,\sqrt5)$ generates the
connected proper subgroup
$H=\{(-u,-v,u,w):u,v,w\in\T\}$, so $H_0=H$ and $m=1$. The lattice vector
$a=(0,1,0,0)\in T_0H_0\cap\Z^4$ gives
\[
m\bigl(\langle a_t,\beta\rangle+\langle a_\omega,\alpha\rangle\bigr)
=\langle(0,1),\beta\rangle=\sqrt5\notin\Z,
\]
contradicting \eqref{eq:arithmetic-rigidity}. Hence no nonzero
$f\in\mathcal S(\R^2)$ yields a linear dependence for this $\Lambda$.
\end{example}

\begin{example}
Let
\[
\Lambda=\bigl\{((0,0),(0,0)),\;((1,0),(0,0)),\;((0,1),(0,0)),\;((0,0),(1,0)),\;(\alpha,\beta)\bigr\}
\subset\R^2\times\R^2,
\]
with $\alpha=(\tfrac12,\sqrt2)$ and $\beta=(0,\sqrt3)$. Then
$\gamma=(-\alpha,\beta)=(-\tfrac12,-\sqrt2,0,\sqrt3)$ has orbit closure
$H=\{0,\tfrac12\}\times\T\times\{0\}\times\T$, so $H_0=\{0\}\times\T\times\{0\}\times\T$
and $m=2$. The vector $a=(0,1,0,0)\in T_0H_0\cap\Z^4$ yields
\[
m\bigl(\langle a_t,\beta\rangle+\langle a_\omega,\alpha\rangle\bigr)
=2\langle(0,1),\beta\rangle=2\sqrt3\notin\Z,
\]
so \eqref{eq:arithmetic-rigidity} fails and no nonzero Schwartz window can
produce a dependence for this configuration.
\end{example}

\begin{remark}
The preceding examples show that in dimension $2$ the infinite proper case of
Theorem~\ref{thm:trichotomy}, the arithmetic
condition from Proposition~\ref{prop:arithmetic} can rule out entire infinite
proper orbit closures outright. In particular, the obstruction is genuinely
stronger in higher dimension than in the one-dimensional examples where the
integrality condition is often automatic.
\end{remark}

\section*{Acknowledgments}
I gratefully acknowledge support from NSF grant DMS-2205852. I also thank Chris Heil, Kasso Okoudjou, and Revati Jadhav for their valuable feedback on earlier versions of this manuscript and for pointing me toward the winding number argument used in the arithmetic rigidity section.

\bibliographystyle{plain}
\bibliography{hrt}

@article{Gabardo2004,
  author    = {Jean-Pierre Gabardo and Deguang Han},
  title     = {The uniqueness of the dual of Weyl-Heisenberg subspace frames},
  journal   = {Appl. Comput. Harmon. Anal.},
  volume    = {17},
  number    = {2},
  pages     = {226--240},
  year      = {2004},
  note      = {MR2082160}
}

@article{Gabardo2012,
  author    = {Jean-Pierre Gabardo},
  title     = {Convolution inequalities in locally compact groups and unitary systems},
  journal   = {Numer. Funct. Anal. Optim.},
  volume    = {33},
  number    = {7-9},
  pages     = {1005--1030},
  year      = {2012},
  note      = {MSC 42C15, MR2966142}
}

@article{Gabardo2014Measure,
  author    = {Xiaoye Fu and Jean-Pierre Gabardo},
  title     = {Measure of self-affine sets and associated densities},
  journal   = {Constr. Approx.},
  volume    = {40},
  number    = {3},
  pages     = {425--446},
  year      = {2014},
  note      = {MSC 28A78}
}

@article{Gabardo2015Gabor,
  author    = {Jean-Pierre Gabardo and Chun-Kit Lai and Yang Wang},
  title     = {Gabor orthonormal bases generated by the unit cubes},
  journal   = {J. Funct. Anal.},
  volume    = {269},
  number    = {5},
  pages     = {1515--1538},
  year      = {2015},
  note      = {MSC 42B05}
}

@article{Gabardo2015SelfAffine,
  author    = {Xiaoye Fu and Jean-Pierre Gabardo},
  title     = {Self-affine scaling sets in $\mathbb{R}^2$},
  journal   = {Mem. Amer. Math. Soc.},
  volume    = {233},
  number    = {1097},
  pages     = {vi+85 pp.},
  year      = {2015},
  isbn      = {978-1-4704-1091-9},
  note      = {MSC 42-02}
}

@article{Gabardo2018Sampling,
  author    = {Jean-Pierre Gabardo},
  title     = {Sampling and interpolation in weighted $L^2$-spaces of band-limited functions},
  journal   = {Sampl. Theory Signal Image Process.},
  volume    = {17},
  number    = {2},
  pages     = {197--224},
  year      = {2018},
  note      = {MSC 42C15}
}

@article{Gabardo2018Weighted,
  author    = {Jean-Pierre Gabardo},
  title     = {Weighted convolution inequalities and Beurling density},
  journal   = {Contemp. Math.},
  volume    = {706},
  publisher = {American Mathematical Society},
  pages     = {175--200},
  year      = {2018},
  isbn      = {978-1-4704-3619-3},
  note      = {MSC 42C15}
}

@article{Gabardo2019Decomposition,
  author    = {Xiaoye Fu and Jean-Pierre Gabardo},
  title     = {Decomposition of integral self-affine multi-tiles},
  journal   = {Math. Nachr.},
  volume    = {292},
  number    = {6},
  pages     = {1304--1314},
  year      = {2019},
  note      = {MSC 28A80}
}

@article{Gabardo2020Frames,
  author    = {Jean-Pierre Gabardo and Deguang Han},
  title     = {Frames and finite-rank integral representations of positive operator-valued measures},
  journal   = {Acta Appl. Math.},
  volume    = {166},
  pages     = {11--27},
  year      = {2020},
  note      = {MSC 42C15}
}

@article{Gabardo2020LocalFourier,
  author    = {Jean-Pierre Gabardo},
  title     = {Local Fourier spaces and weighted Beurling density},
  journal   = {Adv. Oper. Theory},
  volume    = {5},
  number    = {3},
  pages     = {1229--1260},
  year      = {2020},
  note      = {MSC 42C15}
}

@article{Gabardo2020OpenSet,
  author    = {Xiaoye Fu and Jean-Pierre Gabardo and Hua Qiu},
  title     = {Open set condition and pseudo Hausdorff measure of self-affine IFSs},
  journal   = {Nonlinearity},
  volume    = {33},
  number    = {6},
  pages     = {2592--2614},
  year      = {2020},
  note      = {MSC 28A80}
}

@article{Gabardo2021,
  author    = {Jean-Pierre Gabardo and Chun-Kit Lai and Vignon Oussa},
  title     = {On exponential bases and frames with non-linear phase functions and some applications},
  journal   = {J. Fourier Anal. Appl.},
  volume    = {27},
  number    = {2},
  pages     = {Paper No. 9, 23 pp.},
  year      = {2021},
  note      = {MSC 42C15}
}

@article{Gabardo2024,
  author    = {Jean-Pierre Gabardo},
  title     = {The Turán problem and its dual for positive definite functions supported on a ball in $\mathbb{R}^d$},
  journal   = {J. Fourier Anal. Appl.},
  volume    = {30},
  number    = {1},
  pages     = {Paper No. 11, 31 pp.},
  year      = {2024},
  note      = {MSC 43A45, MR4700865}
}

@article{Oussa2024,
  author    = {Vignon Oussa},
  title     = {Orthonormal bases arising from nilpotent actions},
  journal   = {Trans. Amer. Math. Soc.},
  volume    = {377},
  number    = {2},
  pages     = {1141--1181},
  year      = {2024},
  note      = {MR4688545}
}

@article{Oussa2019,
  author    = {Vignon Oussa},
  title     = {Compactly supported bounded frames on Lie groups},
  journal   = {J. Funct. Anal.},
  volume    = {277},
  number    = {6},
  pages     = {1718--1762},
  year      = {2019}
}

@article{Oussa2018,
  author    = {Vignon Oussa},
  title     = {Frames arising from irreducible solvable actions {I}},
  journal   = {J. Funct. Anal.},
  volume    = {274},
  number    = {4},
  pages     = {1202--1254},
  year      = {2018}
}

@incollection{Oussa2017,
  author    = {Vignon S. Oussa},
  title     = {Regular sampling on metabelian nilpotent Lie groups: the multiplicity-free case},
  booktitle = {Appl. Numer. Harmon. Anal.},
  publisher = {Birkh{\"a}user/Springer},
  address   = {Cham},
  year      = {2017},
  pages     = {377--411},
  isbn      = {978-3-319-55549-2, 978-3-319-55550-8}
}

@book{Grochenig2001,
  author       = {Gr{\"o}chenig, Karlheinz},
  title        = {Foundations of Time-Frequency Analysis},
  publisher    = {Birkh{\"a}user},
  address      = {Boston},
  year         = {2001},
  isbn         = {978-0-8176-4022-6}
}

@article{okoudjou2019extension,
  title={Extension and restriction principles for the HRT conjecture},
  author={Okoudjou, Kasso A},
  journal={Journal of Fourier Analysis and Applications},
  volume={25},
  number={4},
  pages={1874--1901},
  year={2019},
  publisher={Springer}
}

@book{oussa2025lieframe,
  author    = {Vignon Oussa},
  title     = {A Bridge Between Lie Theory and Frame Theory: Applications of Lie Theory to Harmonic Analysis},
  publisher = {Wiley},
  year      = {2025},
  month     = {April},
  edition   = {1st}

}

@book{heil2010basis,
  title={A basis theory primer: expanded edition},
  author={Heil, Christopher},
  year={2010},
  publisher={Springer Science \& Business Media}
}

@article{heil1996linear,
  author = {Heil, Christopher and Ramanathan, Jayakumar and Topiwala, Pankaj},
  title = {Linear independence of time-frequency translates},
  journal = {Proceedings of the American Mathematical Society},
  volume = {124},
  number = {9},
  pages = {2787--2795},
  year = {1996},
  publisher = {American Mathematical Society},
  doi = {10.1090/S0002-9939-96-03346-1},
  url = {https://www.ams.org/journals/proc/1996-124-09/S0002-9939-96-03346-1/S0002-9939-96-03346-1.pdf},
}

@article{linnell1999neumann,
  title={Von Neumann algebras and linear independence of translates},
  author={Linnell, Peter},
  journal={Proceedings of the American Mathematical Society},
  volume={127},
  number={11},
  pages={3269--3277},
  year={1999}
}

@article{demeter2012hrt,
  title={Proof of the HRT conjecture for (2, 2) configurations},
  author={Demeter, Ciprian and Zaharescu, Alexandru},
  journal={Journal of Mathematical Analysis and Applications},
  volume={388},
  number={1},
  pages={151--159},
  year={2012},
  publisher={Elsevier},
  doi={10.1016/j.jmaa.2011.11.030}
}

@article{bownik2013linear,
  title={Linear independence of time–frequency translates of functions with faster than exponential decay},
  author={Bownik, Marcin and Speegle, Darrin},
  journal={Bulletin of the London Mathematical Society},
  volume={45},
  number={3},
  pages={554--566},
  year={2013},
  publisher={Oxford University Press},
  doi={10.1112/blms/bds119}
}

@article{okoudjou2025hrt,
  author       = {Kasso A. Okoudjou and Vignon Oussa},
  title        = {The HRT Conjecture for Two Classes of Special Configurations},
  journal      = {Journal of Fourier Analysis and Applications},
  year         = {2025},
  note         = {To appear as a letter to the editor},
  doi          = {10.1007/s00041-025-10181-8},
  url          = {https://arxiv.org/abs/2110.04053v4}
}

@article{benedetto2015linear,
  title={Linear independence of finite Gabor systems determined by behavior at infinity},
  author={Benedetto, John J and Bourouihiya, Abdelkrim},
  journal={The Journal of Geometric Analysis},
  volume={25},
  number={1},
  pages={226--254},
  year={2015},
  publisher={Springer}
}

@book{HofmannMorris2013,
  author    = {Karl H. Hofmann and Sidney A. Morris},
  title     = {The Structure of Compact Groups: A Primer for the Student---A Handbook for the Expert},
  edition   = {3},
  series    = {De Gruyter Studies in Mathematics},
  volume    = {25},
  publisher = {De Gruyter},
  address   = {Berlin/Boston},
  year      = {2013}
}

\end{document}